\baselineskip=15pt plus 2pt
\magnification =1200

 at 10truept
 at 10truept
 at 10truept
 at 10truept

\font\bigbold =cmbx10 scaled\magstep2 

\centerline {\bigbold Kernels and point processes associated with}\par
\centerline {\bigbold Whittaker functions}\par
\vskip.1in
\centerline{ Gordon Blower${}^*$ and Yang Chen}\par
\vskip.05in
\noindent {Department of Mathematics and Statistics, Lancaster University, LA1 4YF, England, UK},\par
 and\par
\noindent {Department of Mathematics, University of Macau, Avenda da Universidade, Taipa, Macau, China}\par
\vskip.05in
\noindent ${}^*$ Corresponding author's email: 
{\tt g.blower@lancaster.ac.uk}\par
\vskip.05in
\centerline{9th October 2015}\par
\vskip.05in
\noindent {\bf Abstract.} This article considers Whittaker's function $W_{\kappa ,\mu }$
where $\kappa$ is real and $\mu$ is real or purely imaginary. Then 
$\varphi (x)=x^{-\mu
-1/2}W_{\kappa ,\mu }(x)$ arises as the scattering function of a continuous time linear
system with state space $L^2(1/2, \infty )$ and input and output spaces ${\bf C}$. The
Hankel operator $\Gamma_\varphi$ on $L^2(0, \infty )$ is expressed as a matrix with
respect to the Laguerre basis and gives the Hankel matrix of moments of a
Jacobi weight $w$. The operation of translating $\varphi$ is equivalent to multiplying $w$
by an exponential factor to give $w_\varepsilon$. The determinant of the Hankel matrix
of moments of $w_\varepsilon$ satisfies the $\sigma$ form of Painlev\'e's
transcendental differential equation $PV$. It is shown that $\Gamma_\varphi$ gives rise to
the Whittaker kernel from random matrix theory, as studied by Borodin and Olshanski (Comm. Math. Phys. 211
(2000), 335--358).\par
\par
\vskip.05in
\noindent {\bf 1. Introduction}\par
\vskip.05in
\noindent The Whittaker function $W_{\kappa ,\mu} $ is the solution of the second order linear differential equation
$$y''+\Bigl( {{-1}\over{4}}+{{\kappa}\over{x}}+{{1/4-\mu^2}\over{x^2}}\Bigr) y=0\eqno(1.1)$$
\noindent that is asymptotic to $x^\kappa e^{-x/2}$ as $x\rightarrow\infty$ through real values, and
 possibly has a logarithmic singularity at $x=0$. The general solution of
(1.1) is given by linear combinations of the Whittaker functions $M_{\kappa, \pm \mu}$. 
See [14, 15, 36 p. 343]
 for basic definitions and properties, such as $W_{\kappa, \mu}=W_{\kappa ,-\mu}$. We consider the case in which $\kappa$ is real,
 and $\mu$ is either real or purely imaginary; 
hence $W_{\kappa ,\mu}(x)$ is real for all $x>0$. 
In random matrix theory [4, 5, 6, 23], kernels such as
$$K(x,y)=\bigl( (\kappa -1/2)^2-\mu^2\bigr)
 {\sqrt{xy}}{{ W_{\kappa -1,\mu}(x)W_{\kappa ,\mu}(y)
-W_{\kappa ,\mu}(x)W_{\kappa -1,\mu}(y)}\over{x-y}} ,\eqno(1.2)$$
\noindent provide integral operators on $L^2((0, \infty ); {\bf C})$, and are
associated with determinantal random point fields.
 The purpose of this note is to provide some transparent proofs of some basic 
properties of these kernels and their associated determinantal random point fields. We
show that the Whittaker kernels are closely related to systems of orthogonal polynomials 
for a Pollaczek--Jacobi type weight.\par 
\indent The classical Jacobi weight function [30] can be translated onto $[0,1]$ to become
$$w(x)=x^b(1-x)^a\qquad (0<x<1),\eqno(1.3)$$
\noindent where $a,b>-1$. We will consider deformations of this weight which arise from multiplying by a
positive factor, for instance $e^{-t/x}$ where $t$ is the deformation parameter.\par 
\vskip.05in
\noindent {\bf Definition} With $a,b>-1$ and $t\geq 0$, we introduce
$$\eqalignno{{\cal D}_N(t;a,b)&=\det\Bigl[ \int_0^1 x^{j+k}x^b(1-x)^a e^{-t/x} \,
dx\Bigr]_{j,k=0}^{N-1}&(1.4)\cr}$$
\noindent or equivalently
$$\eqalignno{{\cal D}_N(t;a,b)&=
{{1}\over{N!}}\int_{[0,1]^N} \prod_{1\leq j<k\leq N} (x_j-x_k)^2\prod_{\ell =1}^N 
x_\ell^b(1-x_\ell )^a e^{-t/x_\ell} \, dx_\ell .&(1.5) \cr}$$
\vskip.05in
\indent In section two, we introduce the notion of a linear system and scattering
function.
In section three, we show how the determinants ${\cal D}_N(t;a,b)$ are related to Whittaker
functions and state the Painlev\'e transcendental differential 
equations [16, 20, 24] that they
satisfy. In section four, we introduce the Whittaker kernels and obtain factorization
theorems in the style of Tracy and Widom [32, 33, 34], which express 
Fredholm determinants of the kernels as determinants of 
Hankel integral operators on $L^2(0, \infty )$. These determinants
do not involve $N$ directly. The case $\kappa =a+1/2$ and $\mu =a$ 
is of particular interest, as in [25], and we consider this in sections 5 and 6.\par
%\vskip.05in
\indent A stochastic point process is a probability measure on the space of point
configurations. See section 5.4 of [11] for the general definition
of finite point processes. The process is said to be determinantal when the correlation
functions are given by Fredholm determinants, as follows.\par
\vskip.05in

 \noindent {\bf Definition} (Determinantal point process) 
Let $S$ be a continuous kernel on $(0, \infty )$ such that\par
\indent  (i) $S(x,y)=\overline{S(y,x)}$ for all $x,y\in (0, \infty )$;\par
\indent (ii) the integral operator with kernel $S$ satisfies 
the operator inequality $0\leq S\leq I$ as self-adjoint operators on $L^2(0,\infty )$;\par
\indent (iii) for all $0<u<v<\infty$, the integral operator with kernel ${\bf I}_{[u,v]} 
S(x,y){\bf I}_{[u,v]}$ is trace class, where ${\bf I}_{[u,v]}$ denotes the indicator
function of $[u,v]$.\par
\noindent Then $S$ gives rise to a determinantal point process on
$(u,v)$ in 
Soshnikov's sense [7, 29]. Let $T$ satisfy $I+T=(I-S)^{-1}$ as operators on $L^2(u,v)$. Then
the probability that there are exactly $N$ points in the
realization, one in each subset $dx_j$ for $j=1,\dots ,N$ and none
elsewhere is equal to
$$\det (I+T)^{-1}\det [T(x_j,x_k)]_{j,k=1}^N dx_1\dots dx_N,\eqno(1.6)$$
\noindent where the $x_j\in [0,1]$.\par
\vskip.05in
\noindent {\bf 2. Linear systems for the Whittaker functions}\par
\vskip.05in

\noindent {\bf Definition} (Linear systems)  Let $H$ be a complex 
separable Hilbert space known as the state space, 
and $H_0$ a finite dimensional complex Hilbert space which serves as the 
input and output space. Let ${\cal L}(H)$ be the space of bounded linear
operators on $H$ with the operator norm. A linear system $(-A,B,C)$ consists of:\par
\indent (i) $-A$, the generator of a $C_0$ (strongly continuous) semigroup $(e^{-tA})_{t\geq 0}$ of bounded linear operators on $H$ such that $\Vert e^{-tA}\Vert\leq Me^{-\omega_0t}$ for all $t\geq 0$ and some $M, \omega_0\geq 0$;\par
\indent (ii) $B:H_0\rightarrow H$ a bounded linear operator;\par
\indent  (iii) $C:H\rightarrow H_0$ a bounded linear operator.\par
\noindent We define the scattering function 
$\varphi\in  C((0, \infty ); {\cal L}(H_0))$ by $\varphi (t)=Ce^{-tA}B$, and the  Hankel operator with
 scattering function $\varphi$ by
$$\Gamma_\varphi f(x)=\int_0^\infty \varphi (x+y)f(y)\, dy\qquad (f\in L^2((0,\infty ); H_0)),
\eqno(2.1)$$
\noindent as in [26]. Note that if the integral 
$\int_0^\infty t\Vert \varphi (t)\Vert^2_{{\cal L}(H_0)}dt$ converges, then 
$\Gamma_{\varphi}$ defines a Hilbert--Schmidt operator. In particular, this holds if $\omega_0>0$.\par 
\indent We also introduce $R_\varepsilon :H\rightarrow H$ by
$$R_\varepsilon =\int_0^\infty e^{-tA}B_\varepsilon 
C_\varepsilon e^{-tA}\, dt,\eqno(2.2)$$
\noindent and this integral plainly converges whenever $\omega_0>0$. See [3] for
some related results.\par
\indent For Whittaker functions, the  basic linear system is the following. 
Let $\varepsilon>0$, $H_0={\bf C}$ and $H=L^2((1/2, \infty ); {\bf C})$ and ${\cal D}(A)=
\{ f(s)\in H: sf(s)\in H\}$. Then
$$\eqalignno{ A:&f(s)\mapsto sf(s),\qquad (f\in {\cal D}(A));&(2.3)\cr
B_\varepsilon:& b\mapsto e^{-\varepsilon s} (s+1/2)^{(\kappa +\mu -1/2)/2}
(s-1/2)^{(-\kappa+\mu -1/2)/2}b, \qquad (b\in {\bf C});\cr
C_\varepsilon: & f(s)\mapsto \int_{1/2}^\infty e^{-\varepsilon s}  
{{(s+1/2)^{(\kappa +\mu -1/2)/2}(s-1/2)^{(-\kappa+\mu -1/2)/2} 
}\over{\Gamma (\mu-\kappa +1/2)}}f(s)ds, \qquad (f\in {\cal D}(A)).
\cr}$$
\noindent {\bf 2.1 Lemma} {\sl (i) Then the scattering function is
$$\varphi_{(\varepsilon )}(x)={{W_{\kappa ,\mu}(x+2\varepsilon )}
\over{(x+2\varepsilon )^{\mu +1/2}}}.\eqno(2.4)$$
\indent (ii) Suppose that $\kappa, \mu\in {\bf R}$ and
 $-\kappa +\mu+1/2>0$. Then $R_\varepsilon$ is
self-adjoint and nonnegative.\par
\indent (iii) Suppose that $\Re \mu >\kappa -1/2$, and either $\varepsilon >0$, or 
$\Re \mu <1/2$ and $\varepsilon =0$. Then $R_\varepsilon$ is trace
class and} 
$$\det (I-\lambda R_\varepsilon )=\det (I-\lambda\Gamma_{\varphi_{(\varepsilon )}})\qquad
(\lambda\in {\bf C}).\eqno(2.5)$$   
\vskip.05in
\noindent {\bf Proof.} (i) This is by direct computation, in which one uses a
 representation formula 
$${{W_{\kappa ,\mu}(x)}\over{x^{\mu +1/2}}}=\int_{1/2}^\infty
 e^{-sx}(s+1/2)^{(\kappa +\mu -1/2)}(s-1/2)^{(-\kappa+\mu -1/2)}{{ds}\over
{\Gamma (-\kappa +\mu +1/2)}}\eqno(2.6)$$
for the Whittaker function from [13, 15].\par
\indent (ii) The operator $R_\varepsilon$ on $L^2(0,\infty )$ is represented by the kernel 
$${{ e^{-\varepsilon s} (s+1/2)^{(\kappa +\mu -1/2)/2}}\over
 {(s-1/2)^{(\kappa-\mu +1/2)/2}}}{{e^{-\varepsilon t} (t+1/2)^{(\kappa +\mu -1/2)/2}}\over{
(t-1/2)^{(\kappa-\mu +1/2)/2}}}{{1}\over{(s+t)\Gamma 
(-\kappa+\mu  +1/2)}}.\eqno(2.7)$$
\noindent The first two factors are multiplication operators, while the
final factor in Carleman's operator $\Gamma$ on $L^2(0, \infty )$ with kernel
$1/(x+y)$, as discussed in [26, p.440]. Now $\Gamma$ is non negative
as an operator, 
hence $R_\varepsilon$ is also non negative by (2.7).
See [17] for more analysis of operators of this form.\par
\indent (iii) One can introduce Hilbert--Schmidt operators 
$\Xi, \Theta :L^2(0, \infty
)\rightarrow L^2(1/2, \infty )$ by $\Theta f=\int_0^\infty
e^{-tA^\dagger} C_\varepsilon^\dagger f(s)\, ds$ and $\Xi f=
\int_0^\infty
e^{-tA} B_\varepsilon f(s)\, ds$ such that $\Gamma_{\varphi_{(\varepsilon
)}}=\Theta^\dagger \Xi$ and $R_\varepsilon=\Xi\Theta^\dagger$; see [3] for more details. The kernel of $\Xi$
is 
$$e^{-\varepsilon s-st}(s+1/2)^{(\kappa +\mu -1/2)/2}
(s-1/2)^{(-\kappa+\mu -1/2)/2}\qquad (s>1/2, t>0),\eqno(2.8)$$
\noindent which is Hilbert--Schmidt, hence $R_\varepsilon$ and likewise 
$\Gamma_{\varphi_{(\varepsilon)}}$ are trace class with
$$\det (I-\lambda R_\varepsilon )=\det (I-\lambda \Xi\Theta^\dagger)=\det (I-\lambda
\Theta^\dagger\Xi )=\det (I-\lambda \Gamma_{\varphi_{(\varepsilon
)}}).\eqno(2.9)$$
\vskip.05in
 \indent For $\alpha >0$, let $L_n^{(\alpha )}(x)$ be the monic generalized Laguerre polynomial of degree $n$, defined by
$$L_n^{(\alpha )}(x)=(-1)^n{{e^{x}}\over{ x^\alpha}} {{d^n}\over{dx^n}} 
\bigl( x^{n+\alpha } e^{-x}\bigr);\eqno(2.10)$$
\noindent then the $n^{th}$ generalized Laguerre function is $\phi_n^{(\alpha )}(x)=e^{-x/2}x^\alpha L^{(\alpha )}_n(x)/\Gamma (n+1)$.\par
\vskip.05in
\noindent {\bf 2.2 Lemma } {\sl Let $\varphi$ be the scattering function $\varphi (x)=W_{\kappa ,\mu}(x)x^{-\mu -1/2}$, and let $w$ be the translated Jacobi weight 
$$w(\xi )=\xi^{\mu-\kappa -1/2}(1-\xi)^{2\alpha -2\mu +1}\qquad (0<\xi <1)
.\eqno(2.11)$$
\noindent Then the operation of $\Gamma_\varphi$ on the generalized Laguerre basis is represented by a matrix of moments for the weight $w$, so}
$$\bigl\langle \Gamma_{\varphi }\phi_\ell^{(\alpha )}, \phi_n^{(\alpha )}\bigr\rangle_{L^2(0, \infty )}={{\Gamma (n+1+\alpha)}\over{\Gamma (n+1)}}
{{\Gamma (\ell+1+\alpha)}\over{\Gamma (\ell +1)}}\int_0^1 
\xi^{\ell +n}w(\xi )d\xi . \eqno(2.12)$$ 
\vskip.05in
\noindent {\bf Proof.} This is suggested by [27]. The Laplace transform of $\phi^{(\alpha )}_n$ satisfies
$$\hat \phi^{(\alpha )}_n(s)={{1}\over{\Gamma (\alpha +n)}}{{(s-1/2)^n}\over{(s+1/2)^{n+1+\alpha}}},
\eqno(2.13)$$
\noindent as one checks by repeatedly integrating by parts. Using the representation formula for $\varphi$ as in (2.6),
 one can express $\langle \Gamma_\varphi \phi_n^{(\alpha )},
 \phi_\ell^{(\alpha )}\rangle $ as an integral with respect to $s$ over $(1/2, \infty )$.
 By changing variables to $\xi =(s-1/2)/(s+1/2)$, one obtains the integral of moments 
with respect to the weight $w$.\par
\vskip.05in
\indent  We now show how the leading minors of the Hankel operator 
$\Gamma_\varphi$ are
related to the Jacobi unitary ensemble. The joint probability density function of the Jacobi
unitary ensemble on $[0,1]^N$ as in [23, 34] is 
$${{1}\over{N!}}{{1}\over{\Gamma (\kappa +\mu +1)^N}}
\prod_{j=0}^{N-1} {{\Gamma (j+1+\alpha )}\over{\Gamma (j+1)}} 
\prod_{0\leq j<k\leq N-1}(x_j-x_k)^2\prod_{j=0}^{N-1}w(x_j).\eqno(2.14)$$
Let $\Delta_N(t)$ be the multiple integral
$$\Delta_N(t)={{1}\over{N!}}{{1}\over{\Gamma (\kappa +\mu +1)^N}}
\prod_{j=0}^{N-1} {{\Gamma (j+1+\alpha )}\over{\Gamma (j+1)}}
\int_{[t,1]^N} 
\prod_{0\leq j<k\leq N-1}(x_j-x_k)^2\prod_{j=0}^{N-1}w(x_j)dx_j,\eqno(2.15)$$
\noindent as in Chen and Zhang [10].
\vskip.05in
\noindent {\bf 2.3 Proposition} {\sl (i) The leading minors of the 
determinant of $\Gamma_\varphi$ satisfy
$$\det \Bigl[ \bigl\langle \Gamma_\varphi \phi_n^{(\alpha )}, \phi_\ell^{(\alpha )}
\bigr\rangle \Bigr]_{\ell ,n=0}^{N-1}=\Delta_N(0).\eqno(2.16)$$
\indent (ii) Let $x_0, \dots , x_{N-1}$ be a sample of $N$ points from the 
Jacobi unitary ensemble. Then the probability of the event  $[x_j\geq t$ for all $j]$ 
equals $\Delta_N(t)/\Delta_N(0).$}\par
\vskip.05in
\noindent {\bf Proof.} (i) This identity follows from Lemma 2.2 and the 
Heine--Andreief identity from [30].\par
\indent (ii) Let $P^{(a,b)}_n(x)$ be the monic Jacobi polynomial of degree $n$ for the weight 
$w(x)=x^b(1-x)^a$ on $[0,1]$, where we choose $b=\mu+\kappa -1/2$ and 
$a=2\alpha -2\mu +1$. See [30, page 58]. Then with the constants 
$\gamma_j=\int P_j^{(a,b)}(x)^2w(x)dx$, the kernel
$$J_N(x,y)=\sum_{j=0}^{N-1} {{P_j^{(a,b)}(x)P_j^{(a,b)}(y)}\over{\gamma_j}}\eqno(2.17)$$
\noindent defines a self-adjoint operator on $L^2(w, [0,1])$ such that $0\leq J_N\leq I$. Hence 
$${{\Delta_N(t)}\over{\Delta_N(0)}}=\det \bigl(I-J_N{\bf I}_{(0,t)}\bigr).\eqno(2.18)$$
\vskip.05in
\noindent {\bf 2.4 Proposition} {\sl Let $$H_N(t)=t(1-t){{d}\over{dt}}\log 
\Delta_N(t).\eqno(2.19)$$ 
Then $\sigma (t)=H_N(t)-d_1-td_2$ satisfies the $\sigma$ form of Painlev\'e's 
transcendental differential equation 
$PVI$, so
$$\sigma'\bigl[ t(t-1)\sigma''\bigr]^2+\bigl[ 2\sigma'(t\sigma'-\sigma
)-(\sigma')^2-\nu_1\nu_2\nu_3\nu_4\bigr]^2=
(\sigma'+\nu_1^2)(\sigma'+\nu_2^2)(\sigma'+\nu_3^2)(\sigma'+\nu_4^2)\eqno(2.20)$$
\noindent where $\nu_1=(a+b)/2$, $\nu_2=(b-a)/2,$ $\nu_3=\nu_4 =(2N+a+b)/2$ with initial
conditions $\sigma (0)=d_2$ and $\sigma'(0)=d_1.$}\par

\vskip.05in
\noindent {\bf Proof.} See [10], and [20].\par
\vskip.05in

\noindent {\bf 3. Determinant formulas for the Pollaczek--Jacobi type weight}\par
\vskip.05in
\indent In this section, we consider the Pollaczek--Jacobi type weights, and
show that translating the scattering function $\varphi$ has the same effect as
deforming the weights through multiplication by $e^{-2t/x}$.   
The $m^{th}$ moment of the Pollaczek--Jacobi weight is defined to be
$$\mu_m(t;a,b)=\int_0^1 x^mx^b(1-x)^a e^{-t/x} \, dx\qquad (m=0, 1, \dots
).\eqno(3.1)$$
\vskip.05in
\noindent {\bf 3.1 Proposition} {\sl The moments satisfy}
$$\mu_m(t;a,b)=\Gamma (a+1)e^{-t/2}t^{(b+m)/2}W_{-(2a+b+m+2)/2,-(b+m+1)/2}(t).\eqno(3.2)$$ 
\vskip.05in
\noindent {\bf Proof.} Making the change of variable $1/x=s+1/2$ in (2.6) we
have
$$\mu_m(t;a,b)=e^{-t/2}\int_{1/2}^\infty e^{-st}(s-1/2)^{-\kappa +\mu
-1/2}(s+1/2)^{\kappa +\mu -1/2}\, ds\eqno(3.3)$$
\noindent with $\kappa =-(2a+b+m+2)/2$ and $\mu =-(b+m+1)/2$.\par  
See [15, 3.4712, page 367]. With the change
of integration variable $x\mapsto 1/\xi$, 
$$\mu_m(t;a,b)=\int_1^\infty \xi^{-2-a-b-m}(\xi -1)^ae^{-t\xi}d\xi.$$
\vskip.05in
\indent The translation operation $\varphi (t)\mapsto 
\varphi_{(\varepsilon )}(t)=\varphi (t+2\varepsilon )$ replaces $w$  by 
$$w_\varepsilon (x)=w(x)\exp\Bigl(-2\varepsilon\Bigl(
{{1}\over{x}}-{{1}\over{2}}\Bigr)\Bigr).\eqno(3.4)$$
\noindent For this problem, we shall be concerned with 
$$\eqalignno{{}&D_N(\varepsilon )&(3.5)\cr
&=\prod_{j=0}^{N-1} {{\Gamma (j+\alpha +1)}\over{\Gamma (j+1)}}[\Gamma (\kappa
-\mu +(1/2))]^{-N} {{1}\over{N!}}\int_{[0,1]^N} \prod_{0\leq j<k\leq N-1}
(x_k-x_j)^2 \prod_{\ell
=1}^N w_\varepsilon (x_\ell )dx_\ell.\cr}$$
Hence a straightforward change of variables gives
$$D_N(\varepsilon )=C_N(\varepsilon ){\cal D}_N (2\varepsilon ;a,b)
\eqno(3.6)$$
\noindent where ${\cal D}_N (t;a,b)$ is defined in (1.4) and 
$$C_N(\varepsilon )=e^{\varepsilon N}[\Gamma (\kappa -\mu +1/2)]^{-N}
\prod_{j=0}^{N-1} {{\Gamma (j+\alpha +1)}\over{\Gamma
(j+1)}}.\eqno(3.7)$$

\noindent As in Lemma 2.2 and Proposition 2.3, we have
$$D_N(\varepsilon )=\det\Bigl[ \bigl\langle 
\Gamma_{\varphi_{(\varepsilon )} }\phi_n^{(\alpha )},
 \phi_\ell^{(\alpha )}\bigr\rangle \Bigr]_{\ell
,n=0}^{N-1}.\eqno(3.8)$$
 Let the quantity $\tilde H_N(t)$ be defined as follows
$$\tilde H_N(t)=t{{d}\over{dt}} \log {\cal
D}_N(t)-N(N+b+a).\eqno(3.9)$$
\vskip.05in
\noindent {\bf 3.2 Theorem} {\sl Then $\tilde H_N(t)$ satisfies
the Jimbo--Miwa--Okamoto $\sigma$ form of Painlev\'e's $PV$ for a
special choice of parameters, so}   
$$\eqalignno{(t\tilde H_N'')^2&=-4t(\tilde H_N')^3+
(\tilde H_N')^2\bigl( 4\tilde H_N +(b+2a+t)^2
+4N(N+a+b) -4a(b+a)\bigr)\cr
&\qquad+2\tilde H'_N \bigl( -(b+2a+t)\tilde H_N-2N a(N+b+a)\bigr) +
(\tilde
H_N)^2.&(3.10)\cr}$$  
\vskip.05in
\noindent {\bf Proof.} This was found in Chen and Dai [9, p. 2161].\par 
\vskip.05in
\noindent {\bf 3.3 Remarks} (i) Chen and Dai [9, Theorem 6.1] also how that $(\tilde
H_N)_{N=1}^\infty$ satisfies a second order nonlinear difference
equation.\par
\indent (ii) ${\cal D}_N(t;a,b)$ is the Wronskian determinant of $$\{
\mu_{2N-1}(t;a,b), \mu_{2N-1}'(t;a,b),\dots
,\mu_{2N-1}^{(N-1)}(t;a,b)\};\eqno(3.11)$$ 
\noindent thus $\mu_{2N-1}(t;a,b)$ determines 
${\cal D}_N(t;a,b)$.\par
\indent (iii) The Painlev\'e V equation has previously appeared in various
ensembles in random matrix theory. Tracy and Widom obtained $PV$ in [34] from the
Laguerre ensemble, and also from the Bessel ensemble from [32] which is associated with
hard edge distributions. In [34], they also considered the ensembles $U_N$ of $N\times N$
complex unitary matrices $U$ with Haar measure, and obtained $PV$ from the
distribution of ${\hbox{trace}}(U)$. Remarkably, this is related to the
uniform measure on the symmetric group $S_N$ of permutations  as
$N\rightarrow\infty$.  Borodin and Olshanski considered the
pseudo-Jacobi ensemble 
and obtained $PV$ from a Fredholm determinant
associated with the Whittaker functions $M_{\kappa ,\mu}(1/x)$ with
$\kappa$ and $\mu$ real as in [6]. Their results have
applications to the infinite dimensional unitary group $U_\infty$.
Lisovyy considered how the determinant for the hypergeometric kernel degenerates to the
determinant for the Whittaker kernel, and realised $PV$ as a limiting
case of $PVI$; see [22, section 10].\par  
\indent The modified Bessel function as in (6.1) satisfies $K_\mu
(z)=(2z/\pi)^{-1/2}W_{0,\mu }(2z)$ as in [28]; so in view of these results, it is
fitting that $PV$ should emerge from the Whittaker kernel.\par

\vskip.05in
\noindent {\bf 4. The matrix Whittaker kernel}\par
\vskip.05in

\noindent Borodin and Olshanski [5, 6] have considered kernels of the form (1.2). We can factorize the kernel in terms of products of 
Hankel operators, by analogy to the results of [31, 34]. In section 5, we see that
the case $\kappa =a+1/2$ and $\mu =a$ is of particular interest.\par
\vskip.05in
\noindent {\bf 4.1 Proposition} {\sl The kernel satisfies} 
$$\eqalignno{\sqrt{zw}&{{W_{\kappa
,\mu}(z)W_{\kappa -1,\mu}(w)-W_{\kappa -1,\mu}(z)W_{\kappa ,\mu}(w)}\over{
(w-z)}}\cr
&=\int_1^\infty
\Bigl(\sqrt{z}W_{\kappa ,\mu}(sz)\sqrt{w}W_{\kappa -1,\mu}(sw)+\sqrt{z}W_{\kappa -1,\mu}(sz)\sqrt{w}W_{\kappa ,\mu}(sw)\Bigr)
{{ds}\over {2s}}.&(4.1)\cr}$$ 
\vskip.05in
\noindent {\bf Proof.} Note that the left-hand side is a continuously differentiable
function of $(z,w)\in (0,\infty )^2$ and that the left-hand side converges to zero as
$z\rightarrow\infty$ or $w\rightarrow\infty$.\par 
 \indent In the following proof we use that matrices
$$I=\left[\matrix{1&0\cr 0&1\cr}\right],\quad J=\left[\matrix{0&-1\cr 1&0\cr}\right],\quad \tilde J=
\left[\matrix{0&1\cr 1&0\cr}\right].\eqno(4.2)$$
\indent Combining the differential equation (1.1) with the recurrence relation
$$z{{d}\over{dx}}W_{\kappa ,\mu}(z)=(\kappa -z/2)W_{\kappa ,\mu}(z)-
\bigl(\mu^2-(\kappa-1/2)^2\bigr)W_{\kappa -1,\mu}(z),\eqno(4.3)$$
\noindent we obtain the matrix differential equation
$$ z{{d}\over{dz}}W =\bigl((1/2) I+\Omega (z)\bigr)W,\eqno(4.4)$$
%$$z{{d}\over{dz}}\left[\matrix{W_{\kappa ,\mu}(z)\cr W_{\kappa -1,\mu}(z)\cr}\right]=
%\left[\matrix{\kappa -z/2&(\kappa -1/2)^2-\mu^2\cr -1&1-\kappa +z/2}\right] 
%\left[\matrix{W_{\kappa ,\mu}(z)\cr W_{\kappa -1,\mu}(z)\cr}\right].\eqno(4.4)$$
%\noindent which we abbreviate to
%$$ z{{d}\over{dz}}Y =\bigl((1/2) I+\Omega (z)\bigr)Y,\eqno(4.5)$$
%\noindent where 
$$W(z)=\left[\matrix{W_{\kappa ,\mu}(z)\cr W_{\kappa -1,\mu}(z)\cr}\right],\quad \Omega (z)=\left[\matrix{\kappa -1/2 -z/2& (\kappa -1/2)^2-\mu^2\cr -1& 1/2-\kappa +z/2 \cr}\right]  
\eqno(4.5)$$
\noindent in which ${\hbox{trace}}(\Omega )=0$. The eigenvalues of $(1/2)I+\Omega
(0)$ are $(1/2)\pm \mu$, so the eigenvalues differ by an integer if and only if
$2\mu$ is an integer; this characterizes the exceptional case for Birkhoff
factorization into canonical form [35].\par
\indent The system (4.4) resembles the system of matrix differential equation
considered by Tracy and Widom [34], although the trace of the coefficient matrix is
non-zero, so we use a variant on their methods as in [1]. We compute
$$\eqalignno{\Bigl( z{{d}\over{dz}}+ w{{d}\over{dw}}\Bigr)
{{\langle JW(z), W(w)\rangle }\over{z-w}}&={{\langle JzW'(z), W(w)\rangle }\over{z-w}}+{{\langle JW(z), wW'(w)\rangle
}\over{z-w}}\cr
&\quad -{{\langle JW(z), W(w)\rangle }\over{z-w}}\cr
&={{\langle J\Omega (z)W(z), W(w)\rangle }\over{z-w}}+{{\langle JW(z), \Omega
(w)W(w)\rangle }\over{z-w}},&(4.6)\cr}$$
\noindent where $J\Omega (z)+\Omega (w)^tJ=-(1/2)(z-w)\tilde J$, so
$$\Bigl( z{{d}\over{dz}}+ w{{d}\over{dw}}\Bigr)
{{\langle JW(z), W(w)\rangle }\over{z-w}}=-(1/2)\langle \tilde JW(z), W(w)\rangle .\eqno(4.7)$$
\noindent For comparison, we have
$$\eqalignno{\Bigl( z{{d}\over{dz}}+ w{{d}\over{dw}}\Bigr)&  
{{1}\over{2}}\int_1^\infty \langle
\tilde JW(sz), W(sw)\rangle {{ds}\over{s}}\cr
&={{1}\over{2}}\int_1^\infty 
\Bigl(\langle
\tilde JzW'(sz), W(sw)\rangle +\langle JW(sz), wW'(sw)\rangle\Bigr) ds\cr
&={{1}\over{2}}\int_1^\infty {{d}\over{ds}}\langle
\tilde JW(sz), W(sw)\rangle ds\cr
&= -{{1}\over{2}}\langle \tilde JW(z), W(w)\rangle .&(4.8)\cr}$$
\noindent Hence the sum
$${{\langle JW(z), W(w)\rangle }\over{z-w}}-
{{1}\over{2}}\int_1^\infty \langle
\tilde JW(sz), W(sw)\rangle {{ds}\over{s}}\eqno(4.9)$$
\noindent is a function of $z/w$, which converges to zero as $z\rightarrow\infty$ or
$w\rightarrow 0$ in any way whatever, so this sum is zero. To obtain the stated result,
we multiply by $\sqrt {zw}$ and rearrange the terms.\par
\vskip.05in 
\indent Let $K$ be the integral operator on
$L^2((a_1,a_2); dx)$ with kernel
$$K(x,y)={{\langle JW(x),W(y)\rangle}\over{x-y}}.\eqno(4.10)$$
\noindent For suitable $a_1$ and $a_2$, we can suppose that $\Vert K\Vert <1$ as an operator on $L^2((a_1,a_2);dx)$ and 
let $S=K(I-K)^{-1}$; so that $(I+S)(I-K)=I$ and $S :L^2((a_1,a_2);dx)\rightarrow
L^2((a_1,a_2);dx)$ has a kernel
$$S(x,y)={{Q(x)P(y)-Q(y)P(x)}\over{x-y}},\eqno(4.11)$$ 
\noindent and we take $S(x,x)=Q'(x)P(x)-Q(x)P'(x)$ on the diagonal.

Let $D(a_1,a_2; K)$ be the Fredholm determinant of $K$, regarded as
a function of the endpoints $a_1$ and $a_2$.\par
\vskip.05in
\noindent {\bf 4.2 Proposition} {\sl (i) The restrictions of the kernels of the operators $S$ and 
$S^2$ to the diagonal satisfy 
$$x{{d}\over{dx}} S(x,x)=-P(x)Q(x)+ (S^2)(x,x).\eqno(4.12)$$ 

\indent (ii) The exterior derivative with respect to the endpoints
satisfies}
$$d\log D(a_1,a_2;K)=S(a_1,a_1)da_1-S(a_2,a_2)da_2.\eqno(4.13)$$  
\vskip.05in
\noindent {\bf Proof.} (i) Let $M$ be the operator of
multiplication by the independent variable $x$ and $D$ the operator of
differentiation with respect to $x$. For a integral operator $T$ let $\delta T=[MD,T]-T$, so that $\delta T$ has kernel
$(x\partial/\partial x+y\partial
/\partial y )T(x,y);$ thus $\delta$ is a pointwise derivation on the kernels,
while $T\mapsto [MD,T]$ is a derivation on operator composition. We use $\doteq$ to
mean that an operator corresponds to a certain kernel. Then 
$$[M,K]\doteq W_{\kappa ,\mu}(x)W_{\kappa -1,\mu}(y)-
W_{\kappa ,\mu}(y)W_{\kappa -1,\mu}(x)\eqno(4.14)$$
\noindent and (4.7) shows that 
$$\delta K\doteq -2^{-1}\langle \tilde J W(x) , W(y)\rangle .$$
\noindent  Then $[M,S](I-K)=(I+S)[M,K]$, so the kernels have
$$[M,S]\doteq Q(x)P(y)-Q(y)P(x)\eqno(4.15)$$
\noindent where $Q=(I-K)^{-1}W_{\kappa ,\mu}$ and $P=(I-K)^{-1}W_{\kappa -1,\mu
}$ are differentiable functions. Hence $S$ is also an
integral operator with kernel of the form (4.11). Also, by some straightforward manipulations, we have 
$$\delta S=(I-K)^{-1}(\delta K)(I-K)^{-1}+S^2,\eqno(4.16)$$
\noindent where
$$(I-K)^{-1}(\delta K)(I-K)^{-1}\doteq 
-(1/2)(P(x)Q(y)+P(y)Q(x))\eqno(4.17)$$
 and a short calculation shows
that $(\delta S)(x,x)=x(d/dx) S(x,x)$. Hence the result.\par 
\indent (ii) This is a standard consequence of the results from [19].\par
\vskip.05in
\indent Let $\phi_{\kappa -1,\mu}(t)=W_{\kappa -1,\mu}(e^t)$ and
 $\phi_{\kappa -1,\mu} (t)=W_{\kappa ,\mu
}(e^t)$ and form the matrix
$$\Phi (t) =\left[\matrix{ 0&0&\phi_{\kappa ,\mu} (t) &0\cr
 0&0&\phi_{\kappa -1,\mu} (t) &0\cr
\phi_{\kappa -1,\mu} (t)&\phi_{\kappa ,\mu} (t)&0&0\cr
0&0&0&0\cr}\right]\eqno(4.18)$$
\noindent and the integral operator on $L^2((0, \infty ); {\bf C})$ with kernel
$$T\doteq 2{{\phi_{\kappa ,\mu} (t)\phi_{\kappa -1,\mu}
(u)-\phi_{\kappa ,\mu} (u)\phi_{\kappa -1,\mu} (t)}\over{e^t-e^u}}.
\eqno(4.19)$$
\vskip.05in
\noindent {\bf 4.3 Corollary} {\sl The Hankel operator with scattering 
function $\Phi$ operates on $L^2((0, \infty ); {\bf C}^4)$ and satisfies}
$$\det (I-T)=\det (I+\Gamma_\Phi ).\eqno(4.20)$$
\vskip.05in
\noindent {\bf Proof.} By Proposition 4.1, with the change of variable
$s=e^u$, we see that the integral operator satisfies 
$T=\Gamma_{\phi_{\kappa
-1,\mu}}\Gamma_{\phi_{\kappa,\mu}}+\Gamma_{\phi_{\kappa
,\mu}}\Gamma_{\phi_{\kappa -1,\mu}}$ where $\Gamma_\phi$ is
in the standard form of a Hankel operator on $L^2(0, \infty)$ as in (2.1). 
By some elementary
determinant manipulations, the Hankel with matrix valued scattering function
satisfies $\det (I+\Gamma_\Phi )=\det
(I-\Gamma_{\phi_{\kappa -1,\mu}}\Gamma_{\phi_{\kappa,\mu}}-
\Gamma_{\phi_{\kappa
,\mu}}\Gamma_{\phi_{\kappa -1,\mu}})$, which produces the stated result.\par
\vskip.05in

\noindent {\bf 4.4 Remark} By [15, 9.237], Whittaker's differential equation generalizes the
differential equation that the associated Laguerre functions satisfy. 
The column vector $Y(z)=\sqrt{z}W(z)$ satisfies $z(d/dz)Y(z)=\Omega (z)Y(z)$, which resembles the differential equation for 
generalized Laguerre functions on [33, page 60]. In Remark 5.2 of [2] we obtained a factorization
theorem for certain Whittaker kernels which expressly excluded the case of generalized Laguerre
functions.\par   
\vskip.05in
\noindent {\bf 5. A special case of the Whittaker kernel}\par

\vskip.05in

\indent For $-1/2<a<1/2$, let $R_a$ be the integral operator on $L^2(0, \infty )$ that has kernel
$$R_a\doteq  {{(y/x)^ae^{-(x+y)/2}}\over{x+y}}\eqno(5.1)$$
\noindent and $R$ be the operator on $L^2((0,\infty ); {\bf C}^2)$ given by the block matrix
$$R=\left[\matrix{0&R_a\cr -R_{-a}&0\cr}\right]\eqno(5.2)$$

\noindent and let $\Psi$ be the matrix-valued scattering function
$$\Psi (t)=\left[\matrix{ 0& \Gamma (1-2a)/(1+t)^{1-2a}\cr -
\Gamma (1+2a)/(1+t)^{1+2a}&0\cr}\right]\qquad (t>0).\eqno(5.3)$$
\vskip.05in

\noindent {\bf 5.1 Proposition} {\sl (i) The operators $R_{\pm a}$ are bounded on $L^2(0,\infty )$ for $\vert a\vert <1/2$.\par 
\indent (ii)  For all $0<a_1<b_1<\infty $, 
the operators $R$ and $\Gamma_\Psi$ are trace class on $L^2((a_1,b_1); {\bf
C}^2)$} and satisfy
$$\det (I-\lambda R)=\det (I-\lambda \Gamma_\Psi )\qquad (\lambda\in {\bf
C}).\eqno(5.4)$$
\indent {\sl (iii) The operator $R$ is skew self-adjoint, and $I+R$ is invertible 
with block matrix form}
$$R(I+R)^{-1}=\left[\matrix{ R_aR_{-a}(I+R_aR_{-a})^{-1}&-R_a(I+R_{-a}R_{a})^{-1}\cr 
-R_{-a}(I+R_aR_{-a})^{-1}& R_{-a}R_a(I+R_{-a}R_a)^{-1}\cr}\right].\eqno(5.5)$$
\indent {\sl (iv) The integral operator $R_aR_{-a}$ is represented by the kernel}
$$\Gamma (1-2a) {{W_{a+1/2,a}(x)W_{a-1/2,a}(y)-W_{a-1/2,a}(x)W_{a+1/2,a}(y)}\over{(x-y)\sqrt{xy}}}.\eqno(5.6)$$

\vskip.05in

\noindent {\bf Proof.} (i) This was proved by Olshanski [25] using a Fourier argument. 
For completeness we give an equivalent proof involving Mellin transforms. First we check that the operators $R_{\pm a}$  are bounded on $L^2(0,\infty )$. The expression
$$\int_0^\infty {{(x/y)^a}\over { (x/y)+1}} {{f(y)dy}\over{y}}={{1}\over{\sqrt{x}}} \int_0^\infty {{(x/y)^a}\over{ \sqrt{(x/y)}+\sqrt{(y/x)}}} {{\sqrt {y} f(y)\, dy}\over{y}}\eqno(5.7)$$
\noindent is a Mellin convolution as in [28], and one shows by a standard contour integration argument that 
$$\int_0^\infty {{z^{a+i\sigma-1/2}dz}\over{z+1}}=\pi{\hbox{sech}}\, \pi (a+i\sigma )\eqno(5.8)$$
\noindent is bounded for all $\sigma \in {\bf R}$ for all $\vert a\vert <1/2$.\par
\indent (ii) This follows as in Lemma 2.1. Letting $M_2$ be the $2\times 2$ complex matrices with Hilbert--Schmidt
norm, we introduce the state space $H=L^2((\omega_0,\infty ); M_2)$ and the input and output space $H_0={\bf C}^2$, then
$$\eqalignno{A: {\cal D}(A)\rightarrow H:& \quad f(x)\mapsto xf(x);\cr
B_{\omega_0} :H_0\rightarrow H:& \quad b\mapsto \left[\matrix{ 0&e^{-x/2}x^{-a}\cr -e^{-x/2}x^a&0\cr}\right] b;\cr
C_{\omega_0}:H\rightarrow H_0:&\quad  g(y)\mapsto \int_{\omega_0}^\infty\left[\matrix{
e^{-y/2}y^{-a}&0\cr 0&  e^{-y/2}y^{a}\cr}\right] g(y)\, dy.&(5.9)\cr}$$
\indent Then $\Vert e^{-tA}\Vert\leq e^{-t\omega_0}$ for all $t>0$, and the scattering function is 
$$C_{\omega_0} e^{-tA}B_{\omega_0}=\int_{\omega_0}^\infty \left[\matrix{ 0& e^{-ty-y}y^{-2a}\cr
 -e^{-ty-y}y^{2a}&0\cr}\right] dy;\eqno(5.10)$$
\noindent and we note that 
$$C_{\omega_0}e^{-tA}B_{\omega_0}\rightarrow \Psi (t)\eqno(5.11)$$
\noindent as $\omega_0\rightarrow 0+$, the corresponding $R$ operator has kernel
$$R=\int_0^\infty e^{-tA}B_{\omega_0}C_{\omega_0}e^{-tA}dt;\eqno(5.12)$$
\noindent so as $\omega_0\rightarrow 0$, we obtain
$$\eqalignno{R&\doteq\int_0^\infty \left[\matrix{ 0& e^{-t(x+y)-(x+y)/2}(y/x)^a\cr -(x/y)^a e^{-t(x+y)-(x+y)/2}&0\cr}\right]\, dt\cr
&\doteq \left[\matrix{0& (y/x)^a e^{-(x+y)/2}/(x+y)\cr -(x/y)^a 
e^{-(x+y)/2}/(x+y)&0\cr}\right].&(5.13)\cr}$$
\indent (iii) Note that $R_a^*=R_{-a}$ hence $R^*=-R$ and $I+R_{-a}R_a$ is invertible. By elementary row operations one checks that 
$$\left[\matrix{ I& R_a\cr -R_{-a}&I\cr}\right]^{-1}=\left[\matrix{ (I+R_aR_{-a})^{-1}&
-R_a(I+R_{-a}R_a)^{-1}\cr R_{-a}(I+R_aR_{-a})^{-1}& (I+R_{-a}R_a)^{-1}\cr}\right],\eqno(5.14)$$
\noindent in which we observe that the indices $a$ and $-a$ alternate.\par
\indent (iv) We have
$$\eqalignno{R_aR_{-a}&\doteq \int_0^\infty {{(x/y)^ae^{-(x+y)/2}}\over{x+y}}{{(y/z)^{-a}e^{-(y+z)/2}}\over{y+z}}dy\cr
&\doteq {{(xz)^ae^{-(x+z)/2}}\over{z-x}}\int_0^\infty 
\Bigl( {{1}\over{x+y}}-{{1}\over{y+z}}\Bigr) y^{-2a}e^{-y}\, dy&(5.15)\cr}$$
\noindent where [14, 14.2 (17) , and page 431] gives the final integral in terms of an incomplete Gamma function which reduces to Whittaker's function
$$\int_0^\infty {{y^{-2a}e^{-y}}\over {x+y}} dy=\Gamma (1-2a)x^{-a-1/2}e^{x/2}W_{a-1/2,a}(x).
\eqno(5.16)$$
We also have the identity from [14, p 432]
$$x^ae^{-x/2}=x^{-1/2}W_{a+1/2,a}(x).\eqno(5.17)$$

\noindent The result follows on substituting (5.16) and (5.17) into (5.15).\par
\vskip.05in

\vskip.05in

\noindent {\bf 6. Diagonalizing the Whittaker kernel}\par
\vskip.05in
\indent The Whittaker's functions are related to Bessel's functions and Laguerre's functions. 
Let $K_\nu (x)$ be the modified Bessel function of the second kind given by 
$$K_\nu (x)=\int_0^\infty \cosh (\nu t) e^{-x\cosh t}\, dt,\eqno(6.1)$$
\noindent which is holomorphic on the open right half plane $\{ x: \Re x>0\}$ and decays rapidly 
as $x\rightarrow\infty$ through real values. When $\nu =im$ is purely imaginary, one refers to
McDonald's function [28]. \par
\indent Let $L_a$ be the differential operator 
$$L_a f=-{{d}\over{dx}}\Bigl(x^2{{df}\over{dx}}\Bigr) +\bigl( x/2+a\bigr)^2  f(x).\eqno(6.2)$$
\noindent Then $L_{-\kappa}$ has eigenfunctions 
$$f_{\kappa ,m}(x)=x^{-1}W_{\kappa, im}(x)\qquad (m\geq 0),\eqno(6.3)$$ 
so that 
$$L_{-\kappa} f_{\kappa ,m}(x)=\bigl( 1/4+ m^2+\kappa^2 \bigr) f_{\kappa ,m}(x).\eqno(6.4)$$
\vskip.05in
\noindent {\bf 6.1 Proposition} {\sl (i) The differential operator $L_{-a}$ commutes with 
$R_aR_{-a}$ on\par
\noindent $C_c^\infty ((0,\infty ); {\bf C})$, and}
$$L_{-a}R_a=R_aL_{a}.\eqno(6.5)$$
\indent{\sl (ii) $R$ can be expressed as a diagonal operator with respect to 
$(f_{\pm a, m})_{m\geq 0}$.}\par
\vskip.05in
\noindent {\bf Proof.} (i) Suppose that $f$ is smooth and has compact support inside $(0, \infty )$. Then we can repeatedly integrate by parts the integral
$$\int_0^\infty {{(y/x)^ae^{-(x+y)/2}}\over{(x+y)}}\Bigl( - {{d}\over{dy}}\Bigl( y^2 {{df}\over{dy}}\Bigr)+ (y/2+a)^2f(y)\Bigr) dy\eqno(6.6)$$
\noindent without introducing boundary terms, to obtain 
$$\eqalignno{\int_0^\infty {{(x/y)^ae^{-(x+y)/2}}\over{x+y}}&\Bigl( {{y^2}\over{4}}+a^2+{{y^2}\over{(x+y)^2}}+ay\cr
&+{{y^2}\over{x+y}}+{{2ay}\over{x+y}}-a-y+{{y^2}\over{(x+y)^2}}-{{2y}\over{x+y}}\Bigr) f(y)\,dy.&(6.7)\cr}$$
\indent  After a little reduction, one shows that this coincides with 
$$\Bigl( -x^2{{d^2}\over{dx^2}}-2x{{d}\over{dx}}+(x/2-a)^2\Bigr) \int_0^\infty {{(y/x)^ae^{-(x+y)/2}}\over{x+y}}f(y)\, dy.\eqno(6.8)$$
\noindent Likewise, we have $L_aR_{-a}=R_{-a}L_{-a}$; so $L_{-a}$ commutes with $R_aR_{-a}$.\par
\indent (ii) Erdelyi [14, 14.3 (53)] quotes the following formula for the Stieltjes transform
$$\int_0^\infty {{y^{-1-a}e^{-y/2}}\over{x+y}}W_{-a,im}(y)dy=\Gamma ((1/2)-a+im)\Gamma
((1/2)-a-im)x^{-a-1}e^{x/2}W_{a,im}(x),\eqno(6.9)$$
\noindent so $f_{-a,m}(x)=x^{-1} W_{-a,im}(x)$ satisfies 
$$R_af_{-a,m}(x)=\Gamma ((1/2)-a+im)\Gamma ((1/2)-a-im)f_{a,m}(x).
\eqno(6.10)$$
\noindent  It follows that $f_{-a,m}$ is an eigenvector for $R_{-a}R_a$. Wimp [37] 
considered a Fourier--Plancherel formula which decomposes $L^2((0, \infty ); {\bf C})$ as a direct integral with respect to $f_{a,m}(x)$, where $f_{a,m}(x)$ 
are eigenfunctions of $L_{-a}$, so there is a transform pair
$$\eqalignno{g(m)&=\Gamma (1/2-\kappa -im)\Gamma (1/2-\kappa +im)\int_0^\infty f(t)f_{\kappa
,m}(t)\, dt,\cr
f(x)&={{1}\over{\pi^2}}\int_0^\infty t\sinh (2\pi t) f_{\kappa ,t}(x) g(t)\,
dt.&(6.11)\cr}$$

Hence we can take $(f_{a,m})_{m\geq 0}$ to be a basis of $L^2((0,\infty ); {\bf C})$, and 
$(f_{-a,m})_{m\geq 0}$ to be a basis of another copy of $L^2((0, \infty ); {\bf C})$ and diagonalize $R$ with respect to the combined basis 
$(f_{\pm a, m})_{m\geq 0}$ of $L^2((0,\infty ); {\bf C}^2)$.
\vskip.05in

\noindent {\bf 6.2 Corollary} {\sl For $-1/2<a<1/2$, the kernel
$$K\doteq {{\Gamma (1-2a)\cos^2\pi a }\over{ \pi^2}}
{{W_{a+1/2,a}(x)W_{a-1/2,a}(y)-W_{a-1/2,a}(x)W_{a+1/2,a}(y)}\over{\sqrt{xy}(x-y)}}
\eqno(6.12)$$
\noindent gives a determinantal point process on $[0,\infty )$ such that random
points $\lambda_j$, ordered by 
 $0\leq \lambda_1\leq \lambda_2\leq \dots $, satisfy} 
$${\hbox{Pr}}[\lambda_1\geq s]= \det (I-K{\bf I}_{[0,s]}).\eqno(6.13)$$ 
\vskip.05in
\noindent {\bf Proof.} Note that 
$$\Gamma (a+1/2+im)
 \Gamma (a+1/2-im)\Gamma (-a+1/2+im)\Gamma
(-a+1/2-im)=2\pi^2/(\cos 2\pi a +\cosh 2\pi m)$$
\noindent so
$$\Vert R_aR_{-a}\Vert_{{\cal L}(L^2(0,\infty ))}
=\pi^2\sec^2\pi a.\eqno(6.14)$$
\noindent Hence we have $0\leq K\leq I$ as operators on $L^2(0, \infty )$,
and $K(x,y)$ has a continuous kernel. Hence by Mercer's theorem, $K$
restricts to a trace class operator on $L^2(0,b)$ for all $0<b<\infty$ with
$$ {\hbox{trace}}\, (K)=\int_0^b K(x,x)\, dx.\eqno(6.15)$$
\noindent Note that the right-hand side is finite for all $b>0$, but diverges to
$\infty$ as $b\rightarrow\infty$, hence $\det (I-K{\bf I}_{[0,b]})\rightarrow 0$
as $b\rightarrow\infty$.\par
\indent By Soshnikov's theorem [29], there is a point process on $[0,b]$
with $K$ as the generating kernel. The point process distributes random $x$ in
$[0,b]$ such that only finitely many $x$ can lie in each Borel subset of $[0,b]$,
and the joint distribution of the points is specified as follows. 
 Let $B_j$ $(j=1, \dots ,m)$ be disjoint Borel subsets of $[0,b]$ and
$n_j=\sharp \{ x\in B_j\}$ the number of random points that lie in $B_j$. Then the
joint probability generating function of the random variables $n_j$ is 
$${\bf E}\prod_{j=1}^m z_j^{n_j}=\det \Bigl(I+{\bf I}_{[0,b]} \sum_{j=1}^m (z_j-1)K{\bf
I}_{B_j}\Bigr)\qquad (z_j\in {\bf C}).\eqno(6.16)$$
\noindent Equivalently, as in [7, p. 597] we compress $K$ to 
${\bf I}_{[0,b]}K{\bf I}_{[0,b]}$ on $L^2({[0,b]})$ and introduce\par
\noindent $T_{[0,b]}=
{\bf I}_{[0,b]}K{\bf I}_{[0,b]}(I-
{\bf I}_{[0,b]}K{\bf I}_{[0,b]})^{-1}$ so that 
$$\det (I+T_{[0,b]})=\bigl(\det (I-{\bf I}_{[0,b]}K{\bf I}_{[0,b]})\bigr)^{-1}.\eqno(6.17)$$ 
\noindent For each integer $\ell\geq 0$, the $\ell$-point correlation function is defined to be the positive
symmetric function 
$$\rho_\ell (x_1, \dots ,x_\ell )=\det (I+T_{[0,b]})^{-1}\det\bigl[
T_B(x_j,x_k)\bigr]_{j,k=1}^\ell\qquad (x_j\in [0,b]) ,\eqno(6.18)$$
\noindent such that 
$${{\hbox{Pr}}}({\hbox{there are exactly}}\quad\ell\quad{\hbox{particles in}}\quad
B)={{1}\over{\ell !}}\int_{B^\ell}\rho_\ell (x_1, \dots ,x_\ell )dx_1\dots
dx_\ell\eqno(6.20)$$ 
\noindent for all Borel subsets $B$ of $(0, b)$.\par
\vskip.05in
\noindent {\bf Remarks 6.3} (i) In particular [28, (7.4.25)] we can write 
$$f_{0,m}(x)={{1}\over{\sqrt{\pi x}}}K_{im}(x/2)={{1}\over{\sqrt{2}}}
\int_1^\infty P_{-1/2+im}(s)e^{-sx/2}ds \eqno(6.21)$$
\noindent in terms of the MacDonald and associated Legendre functions. 
The function $W_{0,im}$ also occurs in the spectral decomposition of the Laplace operator
over the fundamental domain that arises from the action of $SL(2, {\bf Z})$ on the upper half 
plane; see [21, p 318] for a discussion of Maass cusp forms. In [2], we considered the Hankel operators that commute with second
order differential operators, and found the case $L_0$ and $R_0$ as 
Q(vii).\par
\indent (ii) Another case of Q(vii) from [2] is $L_{-\kappa}$ commuting with the
Hankel operator with scattering function $x^{-1}W_{\kappa ,1/2}(x)$ on
$C_c^\infty (0, \infty )$. If $a\neq 0$, then $R_a$ is not a
Hankel operator but an operator of Howland's type as in [17], and $L_a$ does
 not commute with $R_{-a}$.\par
\vskip.05in  
\noindent {\bf Acknowledgements} Gordon Blower acknowledges the hospitality of the University of
New South Wales, where part of this work was carried out. Yang Chen would like to thank
the Science and Technology Development Fund of Macau S.A.R. for awarding the grant
 FDCT 077/2012/A3, and the University of Macau for awarding MYRG2014-00011-FST and
MYRG2014-00004-FST.\par
\vskip.05in
\noindent {\bf References}\par
\vskip.05in

\noindent [1] G. Blower, Integrable operators and the squares of Hankel operators,
 {\sl J. Math. Anal. Appl.} {\bf 340} (2008), 943--953.\par
\noindent [2] G. Blower, Hankel operators that commute with second-order differential operators,
{\sl J. Math. Anal. Appl.} {\bf 342} (2008), 601-614.\par
\noindent [3] G. Blower, Linear systems and determinantal point processes,
 {\sl  J. Math. Anal. Appl.} {\bf 355} (2009), 311--334.\par

\noindent [4] A. Borodin and G. Olshanski, Point processes and the infinite symmetric group, {\sl Math. Res. Lett.}  {\bf 5} (1998), 799--816.\par
\noindent [5] A. Borodin and G. Olshanski, Distributions on partitions, 
point processes and the hypergeometric kernel, {\sl Comm. Math. Phys.} 
{\bf 211} (2000), 335--358.\par
\noindent [6] A. Borodin and G. Olshanski, Infinite random matrices and ergodic
measures, {\sl Commun. Math. Phys.} {\bf 223} (2001), 87--123.\par
\noindent [7] A. Borodin and A. Soshnikov, Janossy densities. I. Determinantal
ensembles, {\sl J. Statist. Phys.} {\bf 113} (2003), 1--16.\par
\noindent [8] M. Chen and Y. Chen, Singular linear 
statistics of the Laguerre unitary ensemble and
Painlev\'e III: Double scaling analysis, {\sl J. Math Phys} {\bf 56} (2015), 063506; 14pp.\par 
\noindent [9] Y. Chen and D. Dai, Painlev\'e V and a Pollaczek--Jacobi type orthogonal polynomials,
{\sl J. Approximation Theory} {\bf 162} (2010), 2149--2167.\par

\noindent [10] Y. Chen and L. Zhang, Painlev\'e VI and the unitary Jacobi ensemble, {\sl Stud. Appl.
Math.} {\bf 125} (2010), 91--112.\par
\noindent [11] D.J. Daley and D. Vere-Jones, {\sl An Introduction to the
Theory of Point Processes}, (Springer-Verlag, New York, 1958).\par
\noindent [12] F.J. Dyson, Correlations between eigenvalues of a random matrix, {\sl Commun. Math.
Phys.} {\bf 19} (1970), 235--250.\par
\noindent [13] A. Erd\'elyi, {\sl Higher Transcendental Functions} Volume 1, (McGraw--Hill, New York,
1954).\par
\noindent [14] A. Erd\'elyi, {\sl Tables of Integral Transforms} Volume II, (McGraw--Hill, New York, 
1954).\par
\noindent [15] I.S. Gradsteyn and I.M. Ryzhik, {\sl Table of Integrals, Series and Products,}
 7th Edn (Elsevier/Academic Press, Amsterdam, 2007).\par
\noindent [16] V.I. Gromak, I. Laine and S. Shimomura, {\sl Painlev\'e Differential Equations in the
Complex Plane,} (Walter de Gruyter, 2002).\par
\noindent [17] J.S. Howland, Spectral theory of operators of Hankel type, {\sl Indiana Univ. Math. J.} {\bf 41} (1992), 409--426.\par
\noindent [18] M. Jimbo, Monodromy problem and the boundary condition for some Painlev\'e equations,
{\sl Publ. RIMS, Kyoto Univ} {\bf 18} (1982), 1137--1161.\par 

\noindent [19] M. Jimbo, T. Miwa and K. Ueno, Monodromy preserving deformation of 
linear ordinary differential equations with rational coefficients. I General theory and the
$\tau$-function, 
{\sl Physica D: Nonlinear phenomena} {\bf 2} (1981), 306--352.\par
\noindent [20] M. Jimbo and T. Miwa, Monodromy preserving deformation of 
linear ordinary differential equations with rational coefficients. II, 
{\sl Physica D: Nonlinear phenomena} {\bf 2} (1981), 407--448.\par
\noindent [21] S. Lang, $SL_2({\bf R})$, (Springer, New York, 1985).\par
\noindent [22] O. Lisovyy, Dyson's constant for the hypergeometric
kernel, pp 243-267 in {\sl New trends in quantum integrable systems},
(World Sci. Publ., Hackensack NJ, 2011).\par 
\noindent [23] M.L. Mehta, {\sl Random Matrices}, Third edition, (Elsevier, San Diego, 2004).\par
\noindent [24] K. Okamoto, On the $\tau$-function of the Painlev\'e equations, {\sl Phys. D: Nonlinear
Phenomena} {\bf 2} (1981), 525--535.\par
\noindent [25]  G. Olshanski, Point processes and the infinite 
symmetric group, Part V: Analysis of the matrix Whittaker kernel, 
ArXiv Math 9810014.\par
\noindent  [26] V.V. Peller, {\sl Hankel Operators on Hilbert Space and Their Applications}, (Springer--Verlag, New York, 2003).\par

\noindent [27]  M. Rosenblum, On the Hilbert matrix I, {\sl Proc. Amer. Math. Soc.} {\bf 9} (1958), 137--140.\par 
\noindent [28] I.N. Sneddon, {\sl The Use of Integral Transforms}, (McGraw--Hill, 1972).\par
\noindent [29] A. Soshnikov, Determinantal random point fields, {\sl Russian Math.
Surveys} {\bf 55} (2000), 923--975.\par
\noindent [30] G. Szeg\"o, {\sl Orthogonal Polynomials}, (American Mathematical Society, New York,
1939).\par

\noindent [31] C.A. Tracy and H. Widom, Level spacing distribution and the Airy kernel, 
{\sl Comm. Math. Phys.} {\bf 159} (1994), 151-174.\par
\noindent [32] C.A. Tracy and H. Widom, Level spacing distribution and the Bessel kernel, 
{\sl Comm. Math. Phys.} {\bf 161} (1994), 289-309.\par
\noindent [33] C.A. Tracy and H. Widom, Fredholm determinants, differential equations and
matrix models, {\sl Comm. Math. Phys.} {\bf 163} (1994), 33-72.\par
\noindent [34] C.A. Tracy and H. Widom, Random unitary matrices, 
permutations, and Painlev\'e, {\sl Comm. Math. Phys.} {\bf 207} (1999),
665-685.\par
\noindent [35] H.L. Turrittin, Reduction of ordinary differential equations to Birkhoff
canonical form, {\sl Trans. Amer. Math. Soc.} {\bf 107} (1963), 485-507.\par
\noindent [36] E.T. Whittaker and G.N. Watson, {\sl A course of modern analysis} fourth edition, (Cambridge
University Press, 1996).\par

\noindent [37] J. Wimp, A class of integral transforms, {\sl Proc. Edin. Math. Soc.} {\bf 14} (1964), 33-40.\par
\vskip.1in

\vfill
\eject
\end

\noindent Let $Z(x)$ have Mellin transform $Z^*(s)=\int_0^\infty x^{s-1}
Z(x)\, dx$, and let $f (t)=Z(e^t)$. Then under suitable integrability
conditions stated in [28, p. 273], 
$$f(t)={{1}\over{2\pi i}}\int_{c-i\infty}^{c+i\infty} e^{-st}
 Z^*(s)\, ds.\eqno(4.21)$$
\indent In particular by [15, 6.42(17)], with $Z(x)=e^{-x/2}W_{\kappa ,\mu }(x)$ for $\vert
\Re \mu \vert <1/2$ we have an absolutely convergent integral along the
imaginary axis
$$f(t)={{1}\over{2i}}\int_{-i\infty}^{+i\infty} e^{-st} {{\Gamma
(s+\mu +1/2)}\over{\Gamma (-s+\mu +1/2)}}{{1}\over
{\Gamma (s-\kappa +1)}}\sec \pi (s-\mu )\, ds\qquad (t\in {\bf R})
\eqno(4.22)$$
\noindent where for $\mu$ real and $s$ on the imaginary axis, the exponential
is unimodular, the quotient of Gamma
functions is unimodular, and the final two factors are of exponential decay as
$s\rightarrow\pm i\infty$.\par  
\vskip.05in

\noindent [1] M. Adler and P. van Moerbeke, Matrix integrals, Toda symmetries, Virasoro constraints and orthogonal polynomials, {\sl Duke J. Math.} {\bf 80}
 (1995), 863--891.\par
\noindent [2] Z.D. Bai, Y. Chen and Y.-C. Liang, {\sl Random matrix theory and its applications}
Lecture Notes Series. Institute for Mathematical Sciences. National University of Singapore, 18
(2009).\par
\noindent [3] E.L. Basor, Distribution functions for random variables for ensembles of positive
Hermitian matrices, {\sl Comm. Math. Phys.} {\bf 188} (1997), 327--350.\par
\noindent [4] E.L. Basor and Y. Chen, Toeplitz determinants from compatibility conditions, {\sl
Ramanujan J.} {\bf 16} (2008), 25--40.\par
\noindent [5] E.L. Basor and Y. Chen, Perturbed Hankel determinants, {\sl J. Phys A: Math and
General} {\bf 38} (2005), 10101-10106.\par
\noindent [6] E.L. Basor, Y. Chen and H. Widom, Determinants of Hankel matrices, {\sl J. Funct.
Anal.} {\bf 179} (2001), 214--234.\par

\noindent [7] M. Bertola, B. Eynard and J. Harnad, Semiclassical orthogonal polynomials. matrix models and isomonodromic tau functions, 
{\sl Comm. Math. Phys.} {\bf 263} (2006), 401--437.\par
\noindent [10] L. Boelen, {\sl Discrete Painlev\'e equations and orthogonal polynomials} PhD Thesis,
K.U. Leuven (2010).\par
\noindent [11] L. Boelen, and W.V. Assche, Discrete Painlev\'e equations for recurrence coefficients f
semiclassical Laguerre polynomials, {\sl Proc. Amer. Math. Soc.} {\bf 138} (2010),
1317--1331.\par

\noindent [12] L. Boelen, G. Filipuk and W.V. Assche, Recurrence coefficients of 
generalized Meixner
polynomials and Painlev\'e equations 
{\sl J. Phys. A Math Theor} {\bf 44} (2011), 035202\par
\noindent [19] Y. Chen and M.V. Feigin, Painlev\'e IV and degenerate Gaussian Unitary ensembles, {\sl
J. Phys. A: Math. Gener} {\bf 39} (2006), 12381--12393.\par 
\noindent [20] Y. Chen and J. Griffin, Non linear difference equations arising from a deformation of
the $q$-Laguerre weight, {\sl Indagationes Mathematicae} {\bf 26} (2015), 266--279.\par
\noindent [21] Y. Chen, N.S. Haq and M.R. McKay, Random matrix models, double time Painlev\'e
equations and wireless relaying, {\sl J. Math. Phys.} {\bf 54} (2013), 063506.\par
\noindent [22] Y. Chen, and M.E.H. Ismail, Thermodynamic relations of the Hermitian matrix ensembles,
{\sl J. Phys. A: Math. General} {\bf 30} (1997), 6633--6654.\par
\noindent [23] Y. Chen, and M.E.H. Ismail, Ladder operators and differential equations for orthogonal
polynomials, {\sl J. Phys. A: Math General} {\bf 30} (1997), 7817-7829\par
\noindent [24] Y. Chen, and M.E.H. Ismail, Ladder operators for $q$-orthogonal polynomials, 
{\sl J. Math. Anal. Appl.} {\bf 345} (2008), 1--10.\par
\noindent [25] Y. Chen, and M. Ismail, Jacobi polynomials from compatibility conditions, 
{\sl Proc. Amer. Math Soc.} {\bf 133} (2005), 465--472.\par
\noindent [26] Y. Chen and A. Its, Painlev\'e III and a singular linear statistics in Hermitian
random matrix ensembles, I {\sl J. Approx. Theory} {\bf 162} (2010), 270-297.\par
\noindent [27] Y. Chen and N. Lawrence, On the linear statistics of Hermitian random matrices, {\sl 
J. Phys. A: Math. Gen.} {\bf 31} (1998), 1141--1152.\par
\noindent [28] Y. Chen  and M.R. McKay, Coulomb fluid, Painlev\'e transcendents and the information
theory of MIMO systems, {\sl IEEE Trans. Inf. Theo.} {\bf 58} (2012), 4594--4634.\par
\noindent [29] Y. Chen and G. Pruessner, Orthogonal polynomials with discontinuous weights, {\sl J.
Phys. A: Math. Gene.} {\bf 38} (2005), L191--L198.\par
\noindent [34] P.J. Forrester, The spectrum edge of random matrix ensembles, {\sl Nucl. Phys. B.} {\bf
402} (1993), 709--728.\par
\noindent [35] P.J. Forrester and N.S. Witte, Boundary conditions associated with the Painlev\'e III'
and V evaluations of some random matrix ensembles, {\sl J. Phys. A: Math. Gene.} {\bf 39} (2006),
8983-8995.\par
\noindent [36] P.J. Forrester and C.M. Ormerod, Differential equations for deformed Laguerre polynomials,
{\sl J. Approx. Theo.} {\bf 162} (2010), 653--677.\par
\noindent [40] M.E.H. Ismail and Z.S.I. Mansour, $q$-analogues of Freud weights and nonlinear
difference equations, {\sl Adva. Appl. Math.} {\bf 45} (2010), 518--547.\par
\noindent [41] M.E.H. Ismail, I. Nikolova and P. Simeonov, Difference equations and discriminants for
discrete orthogonal polynomials, {\sl Ramanujan J.} {\bf 8} (2005), 475--502.\par
\noindent [42] A.R. Its, A.G. Izergin, V.E. Korepin and N.A, Slavnov, Difference equations for quantum
correlation functions, {\sl Inte. J. Mode. Phys. B} {\bf 4} (1990), 1003-1037.\par

\noindent [46] K. Johansson, On fluctuations of eigenvalues of random Hermitian matrices, {\sl Duke
Math. J} {\bf 91} (1998), 151--204.\par
\noindent [47] A.B. Kuijlaars and L. Zhang, Singular values of products of Ginibre random matrices,
multiple orthogonal polynomials and hard edge scaling limits, {\sl Commun. Math. Phys.} {\bf 322}
(2014), 759--781.\par
\noindent [49] T. Nagao and M. Wadati, Correlation functions and random matrix ensembles related to
classical orthogonal polynomials, {\sl J. Phys. Soci. Japan} {\bf 60} (1991), 3298--3322.\par
\noindent [50] J.M. Normand, Calculation of some determinants using the $s$-shifted 
factorial, {\sl J.
Phys. A: Math. Gene} {\bf 37} (2004), 5737--5762.\par 
\noindent [52] Y. Ohyama, H. Kawamuko, H. Sakai and K. Okamoto, Studies on the Painlev\'e equations V.
Third Painlev\'e equations of special type $P_{III}(D_7)$ and $P_{III}(D_8)$, {\sl J. Math. Sci.
Univ. Tokyo} {\bf 13} (2006), 145--204.\par
\noindent [55] A. Ronveaux and M.A. Felix, {\sl Heun's differential equations} (Oxford University
Press, 1995).\par

\noindent [60] C. Texier and S.N. Majumdar, Wigner time-delay distribution in chaotic cavities and
freezing transition, {\sl Phys. Rev Lett.} {\bf 110} (2013), 250-602.\par

\noindent which reduce by functional identities to
$$=\sqrt{xy}\Bigl( {{xW'_{\kappa ,\mu}(x)W_{\kappa ,\mu}(y)-W_{\kappa ,\mu}(x)yW'_{\kappa ,\mu}(y)}\over{x-y}}\Bigr)+{{\sqrt{xy}}\over{2}}W_{\kappa ,\mu}(x)W_{\kappa ,\mu}(y).
\eqno(3.19)$$
\noindent Note the close resemblance with $R_{a}R_{-a}$. The following result has the character of several results from [18], in which an integrable kernel is expressed as a product of Hankel operators. \par
\vskip.05in

\noindent {\bf 3.2 Proposition} {\sl The Whittaker kernel satisfies}
$$ \eqalignno{{}&{{e^tW'_{\kappa, \mu}(e^t)
W_{\kappa, \mu}(e^s)-W_{\kappa, \mu}(e^t)e^sW'_{\kappa, \mu}(e^s)}\over{ 2\sinh ((t-s)/2)}}&(3.20)\cr
&=\int_0^\infty 4^{-1}e^{(t+u)/2}W_{\kappa, \mu}(e^{t+u})e^{(s+u)/2}W_{\kappa, \mu}(e^{s+u})\cr
&-\kappa e^{-(t+u)/2}W_{\kappa, \mu}(e^{t+u})
e^{-(s+u)/2}W_{\kappa, \mu}(e^{s+u})+e^{(t+u)/2}W_{\kappa, \mu}(e^{t+u})e^{(s+u)/2}W_{\kappa, \mu}(e^{s+u})du.\cr}$$
\vskip.05in
\noindent {\bf Proof.} Let $f(t)=W_{\kappa,\mu} (e^t)$, so that 
$$f''(t)=f'(t)+\Bigl({{ e^{2t}}\over{4}}-\mu e^t +\kappa \Bigr) f(t).\eqno(3.21)$$
\noindent It follows that
$$\Bigl( {{\partial}\over{\partial t}}+{{\partial}\over{\partial s}}\Bigr)\Bigl\{ e^{-(s+t)/2} {{f(t)f'(s)-f'(t)f(s)}\over{\sinh (t-s)/2}}\Bigr\}$$
$$=e^{t/2}f(t)e^{s/2}f(s)\cosh (s-t)/2-2\mu f(t)f(s).\eqno(3.22)$$  
\noindent hence
$$e^{-(s+t)/2} {{f(t)f'(s)-f'(t)f(s)}\over{\sinh (t-s)/2}}$$
$$=\int_0^\infty \mu f(t+u)f(s+u)du +\int_0^\infty \bigl(f(t+u)f(s+u)e^{s+u} +e^{t+u}f(t+u)f(s+u)\bigr) du\eqno(3.23)$$
\vskip.05in

\noindent Then introducing

\noindent Then the kernel
$${{\sqrt{zw}\langle JY(z), Y(w)\rangle}\over{z-w}}={{W_{\kappa ,\mu}(z)wW_{\kappa -1,\mu}(w)-zW_{\kappa-1 ,\mu}(z)W_{\kappa ,\mu }(w)}\over{
\sqrt{zw}(z-w)}}$$
\noindent satisfies
 $$\eqalignno{\Bigl( z{{\partial}\over{\partial z}}&+w{{\partial }\over{\partial w}}\Bigr){{\sqrt{zw}\langle JY(z), Y(w)\rangle}\over{z-w}}\cr
&= {{1}\over{\sqrt{zw}}}\Bigl\langle \left[\matrix{ 1&-1/2\cr -1/2& {{\mu^2-(\kappa -1/2)^2}\over{zw}}\cr}\right] \left[\matrix{ W_{\kappa ,\mu }(z)\cr 
zW_{\kappa -1,\mu}(z)\cr}\right],\left[\matrix{ W_{\kappa ,\mu }(w)\cr 
wW_{\kappa -1,\mu}(w)\cr}\right]\Bigr\rangle.&(3.19)\cr}$$
\noindent Hence we have 
$${{\sqrt{zw}\langle JY(z), Y(w)\rangle}\over{z-w}}=\int_1^\infty F(zu)^tG(wu)\, du$$
\noindent where $F(u)$ and $G(u)$ are complex vector-valued functions which tend to $0$ as $u\rightarrow\infty$. The matrix
 $$\left[\matrix{ 1&-1/2\cr -1/2& 0\cr}\right]$$
\noindent has eigenvalues $1/2\pm 1/\sqrt{2}$.\par
\indent  Note that $\mu^2-(\kappa -1/2)^2<0$ when $\mu =im$ is purely imaginary valued, so for $z,w\rightarrow 0$, a significant contribution to the kernel arises from
$$(m^2+(\kappa -1/2)^2)\int_1^\infty {{W_{\kappa -1,im}(uz)W_{\kappa -1,im}(uw)}\over{u\sqrt {zw}}}du.$$

\indent Let
$$Z(z)=\left[\matrix{ z^{-1}W_{\kappa ,\mu }(z)\cr 
W_{\kappa -1,\mu}(z)\cr}\right],\quad \Phi(z)=\left[\matrix{\kappa z -z-z^2/2 &
(\kappa -1/2)^2-\mu^2)\cr 
-z^2& z-\kappa z+z^2/2\cr}\right].\eqno(3.15)$$
\noindent Then $z^2Z'(z)=\Phi (z)Z(z)$ and $\Phi (z)$ is a polynomial 
matrix with
${\hbox{trace}}\, \Phi (z)=0$. This gives the standard form of a 
Tracy--Widom system.

\noindent so
$$v(x)=-\log w(x)=-b\log x-a\log (1-x)\eqno(3.8)$$
\noindent is replaced by 
$$v_\varepsilon (x)=-\log w_\varepsilon (x)
=-b\log x-a\log (1-x)+\varepsilon {{1+x}\over{1-x}}.\eqno(3.9)$$
\noindent The modified potential $v_\varepsilon$ is convex on $(0,1)$ since
$${{v'_\varepsilon (x)-v'_\varepsilon (y)}\over{x-y}}
={{b}\over{xy}}+{{a}\over{(1-x)(1-y)}}+
{{2\varepsilon (2-x-y)}\over{(1-x)^2(1-y)^2}}>0,\eqno(3.10)$$
\noindent and also rational in $x$  and $y$.\par

and 
$$C_N(\varepsilon )=e^{\varepsilon M} [\Gamma (\kappa -\mu +(1/2)]^{-N} \prod_{j=0}^{N-1}
{{\Gamma (j+\alpha +1)}\over{\Gamma (j+1)}}.\eqno(3.14)$$